# Mixed Hodge Structures on Generalized Theta Divisors and Graph Motives

Mohammad Reza Rahmati

December 3, 2025

**Abstract**

We study the mixed Hodge structures on the complements of generalized theta divisors inside generalized Jacobians and their natural compactifications. For a smooth or nodal curve $C$ equipped with an effective modulus $m$, the generalized Jacobian $J_m$ is a semi-abelian extension of an abelian variety by a complex torus, and its generalized theta divisor $\Theta_m$ admits a determinantal description originating in the derived pushforward of the Poincaré bundle. By choosing a smooth compactification $\overline{J}_m$ with simple normal crossing boundary and applying Deligne's logarithmic mixed Hodge theory, we obtain explicit weight filtrations on $H^\bullet(J_m \setminus \Theta_m, \mathbb{Q})$.

Our main results express each graded piece $\mathrm{Gr}_k^W H^n(J_m \setminus \Theta_m, \mathbb{Q})$ in terms of the cohomology of strata indexed by the dual graph of $C$, together with combinatorial data coming from the modulus. This yields structural decompositions into contributions from Jacobians of partial normalizations, tori, and finite sets. We compute several classes of examples explicitly, including the nodal cubic and low-genus nodal curves, recovering mixed Tate structures in simple cases and detecting higher complexity in more singular settings.

Motivated by the analogy between the weight spectral sequence for $J_m \setminus \Theta_m$ and the description of graph motives associated to Feynman integrals, we formulate a conjectural correspondence between the graded pieces of $H^\bullet(J_m \setminus \Theta_m)$ and motives associated to decorated dual graphs of $(C, m)$. These results establish a bridge between the geometry of generalized Jacobians, degenerations of theta divisors, and the motivic theory of graph hypersurfaces.

## 1 Introduction

Theta divisors lie at the heart of the geometry of curves, abelian varieties, and their degenerations. For a smooth curve $C$ of genus $g$, the classical theta divisor on the Jacobian $J(C)$ encodes the geometry of special divisors and plays a central role in the study of linear series, Prym varieties, and the Schottky problem. When $C$ acquires singularities or is equipped with a modulus, the associated generalized Jacobian $J_m$ becomes a semi-abelian variety in the sense of Rosenlicht [17], and the corresponding generalized theta divisor $\Theta_m$ reflects the interaction between the abelian and toric parts of $J_m$. Despite the ubiquity of generalized Jacobians in arithmetic geometry, geometric representation theory, and the theory of degenerations, comparatively little is known about the *mixed Hodge structures* arising from the complements $J_m \setminus \Theta_m$.

## Background and motivation

For the classical Jacobian $J(C)$, the complement $J(C) \setminus \Theta$ is well understood: it is an affine variety with pure Hodge structure coming from the torus fibration over $\mathrm{Pic}^{g-1}(C)$. In contrast, for generalized Jacobians—extensions of $J(C)$ by a torus determined by the singularities or the modulus—the cohomology of $J_m \setminus \Theta_m$ is typically *mixed*, and its weight filtration encodes subtle information about degenerations of $C$, monodromy, and the stratification of compactified Jacobians. The construction of $\Theta_m$ itself is highly nontrivial: it arises as the determinantal degeneracy locus of the derived pushforward of the Poincaré bundle [16, 5, 13], and its singularities reflect both the geometry of $C$ and the combinatorics of the modulus.

The mixed Hodge theory of open varieties has been systematically developed by Deligne [8, 9] and refined by Saito's theory of mixed Hodge modules [18, 19]. When applied to a smooth compactification $\overline{J}_m$ with simple normal crossing boundary, this theory produces a weight spectral sequence whose $E_1$-page is expressed in terms of the boundary strata. Compactified Jacobians of nodal curves admit stratifications indexed by the dual graph of $C$ and by the multidegree of torsion-free rank-one sheaves [1, 5, 7, 10], suggesting a close combinatorial relationship between the geometry of $(C, m)$ and the structure of $H^\bullet(J_m \setminus \Theta_m)$.

## Overview of results

In this work we carry out a systematic study of the mixed Hodge structures on $H^\bullet(J_m \setminus \Theta_m)$ for curves with modulus and their degenerations. Our approach is based on the logarithmic de Rham complex of a smooth compactification $\overline{J}_m$, the determinantal description of $\Theta_m$, and the boundary stratification inherited from the dual graph of the curve.

**(1) Determinantal description.** We begin by reviewing the derived pushforward of the Poincaré bundle and the universal resolution

$$\mathbf{R}\pi_*\mathcal{P} \simeq [\mathcal{E}^0 \xrightarrow{\Phi} \mathcal{E}^1],$$

whose degeneracy locus $\det(\Phi) = 0$ recovers the generalized theta divisor (Proposition 2.4). This underlies our later analysis by allowing us to treat $\Theta_m$ as part of a normal crossings divisor after resolving singularities.

**(2) Mixed Hodge structures on $J_m \setminus \Theta_m$.** For a smooth compactification $(\overline{J}_m, D)$, where $D = \overline{\Theta}_m \cup \partial\overline{J}_m$ is a simple normal crossing divisor, we express

$$H^n(J_m \setminus \Theta_m, \mathbb{Q}) = \mathbb{H}^n\left(\overline{J}_m, \Omega^\bullet_{\overline{J}_m}(\log D)\right)$$

and describe all associated graded pieces through Deligne's weight spectral sequence (Corollary 3.4).

**(3) Dual graph expression for weight–graded pieces.** Using the stratification of compactified Jacobians by multidegree and singularity type, we prove a structural theorem (Theorem 5.1) expressing $\mathrm{Gr}^W_k H^n(J_m \setminus \Theta_m)$ as a direct sum of cohomology groups of strata indexed by subgraphs of the dual graph of $C$, together with decorations coming from the modulus. Each stratum is a product of a Jacobian of a partial normalization with a torus.

**(4) Explicit computations.** We compute several examples in detail, including Jacobians of nodal cubics and low-genus singular curves. In simple cases (e.g. the nodal cubic) the cohomology of $J_m \setminus \Theta_m$ is mixed Tate, while for more complicated degenerations it exhibits nontrivial weight interactions.

**(5) Relation to graph motives.** The structure of the weight spectral sequence closely resembles the description of the mixed Hodge structures of graph hypersurface complements developed by Bloch–Esnault–Kreimer [2, 3]. Motivated by this analogy, we construct decorated dual graphs associated to $(C, m)$ and formulate a conjectural comparison between $\mathrm{Gr}^W_\bullet H^\bullet(J_m \setminus \Theta_m)$ and motives of graph hypersurfaces (Conjecture 6.1). This provides a conceptual bridge between the geometry of degenerations of curves and the motivic theory of Feynman integrals.

### Applications

Our analysis has several potential applications:

- In the theory of degenerations of abelian varieties, the mixed Hodge structure of $J_m \setminus \Theta_m$ gives precise information about vanishing cycles and monodromy for families of curves.
- In mathematical physics, the connection with graph motives suggests that certain theta–divisor complements may appear naturally in the study of amplitudes and periods of Feynman integrals.
- In the study of compactified Jacobians and moduli of sheaves, the structural decomposition of $\mathrm{Gr}^W_\bullet$ aligns with wall-crossing behavior and with the geometry of semi-abelian pairs.
- In arithmetic geometry, these mixed Hodge structures provide new examples of geometric objects whose motivic realizations interpolate between classical abelian varieties and graph hypersurfaces.

Altogether, this work aims to establish a unified framework connecting the geometry of generalized theta divisors, the combinatorics of dual graphs, and the motivic structures governing Feynman graph hypersurfaces.

## 2 Generalized Jacobians and theta divisors

Let $C$ be a smooth projective connected curve of genus $g \geq 1$ over $\mathbb{C}$, and let $m = \sum_{x \in S} m_x[x]$ be an effective divisor (the *modulus*) with support $S \subset C$. Denote by $J(C)$ the Jacobian of $C$ and by $J_m$ the generalized Jacobian of $(C, m)$ in the sense of Rosenlicht [17, 11].

**Definition 2.1.** The *generalized Jacobian* $J_m$ is the commutative algebraic group representing the functor which assigns to a $\mathbb{C}$–scheme $T$ the group of degree zero divisors on $C \times T$ modulo rational equivalence by principal divisors $\mathrm{div}(f)$ such that

$$f \equiv 1 \mod m, \tag{2.1}$$

i.e. $f - 1$ vanishes along $m$.

Rosenlicht showed that $J_m$ fits into an extension of algebraic groups

$$0 \longrightarrow G_m \longrightarrow J_m \longrightarrow J(C) \longrightarrow 0, \tag{2.2}$$

where $G_m$ is a connected commutative affine algebraic group (a product of copies of $\mathbb{G}_a$ and $\mathbb{G}_m$) determined by $(C, m)$ [17, 15, 12].

Fix a base point $P_0 \in C \setminus S$. For $d \in \mathbb{Z}$, let $\mathrm{Div}^d(C)$ be the space of degree-$d$ divisors, and $\mathrm{Pic}^d(C)$ the Picard variety. Let $\mathrm{Pic}^d_m(C)$ denote the generalized Picard variety with modulus $m$ (cf. [17]).

**Definition 2.2.** The *generalized Abel map with modulus $m$* is

$$\alpha_{m,d} \colon C^{(d)} \longrightarrow \mathrm{Pic}^d_m(C), \qquad D \longmapsto \mathcal{O}_C(D - dP_0), \tag{2.3}$$

where $C^{(d)}$ denotes the $d$-th symmetric power of $C$, and the modulus condition (2.1) is imposed on rational equivalence.

For $d = g-1$ and a suitable choice of base point we recover a generalized theta divisor.

**Definition 2.3.** The *generalized theta divisor* $\Theta_m \subset \mathrm{Pic}^{g-1}_m(C)$ is the Zariski closure of the image

$$\Theta_m := \overline{\alpha_{m,g-1}(C^{(g-1)} \setminus \mathrm{Supp}(m))}. \tag{2.4}$$

In the case $m = 0$ and degree $g-1$ this recovers the usual theta divisor on $J(C)$.

We now recall the determinantal description of $\Theta_m$. Fix a line bundle $L \in \mathrm{Pic}^{g-1}_m(C)$ and consider the complex

$$R\Gamma(C,L) \colon H^0(C,L) \xrightarrow{\delta_L} H^1(C,L), \tag{2.5}$$

in degrees 0 and 1. The theta divisor appears as the degeneracy locus of a universal morphism obtained by globalizing (2.5).

Let $\pi \colon C \times \mathrm{Pic}^{g-1}_m(C) \to \mathrm{Pic}^{g-1}_m(C)$ be the projection, and let $\mathcal{P}$ be a Poincaré line bundle on $C \times \mathrm{Pic}^{g-1}_m(C)$ normalized along $P_0 \times \mathrm{Pic}^{g-1}_m(C)$.

Define the derived pushforward

$$\mathbf{R}\pi_*\mathcal{P} \simeq [\mathcal{E}^0 \xrightarrow{\Phi} \mathcal{E}^1] \tag{2.6}$$

in the derived category $D^b(\mathrm{Coh}(\mathrm{Pic}^{g-1}_m(C)))$, where $\mathcal{E}^i$ are vector bundles and $\Phi$ is a morphism of vector bundles representing $R\Gamma(C,L)$ fiberwise, after possibly twisting $\mathcal{P}$ by a fixed line bundle to kill higher cohomology. See [16, Ch. III] and [13, 5] for details in the compactified case.

**Proposition 2.4** (Determinantal description of $\Theta_m$)**.** *There exists a representation* (2.6) *such that the theta divisor $\Theta_m$ is the degeneracy locus*

$$\Theta_m = \{L \in \mathrm{Pic}^{g-1}_m(C) \mid \det(\Phi_L) = 0\} = V(\det(\Phi)), \tag{2.7}$$

*where* $\det(\Phi)$ *is a global section of* $\det(\mathcal{E}^1) \otimes \det(\mathcal{E}^0)^\vee$. *In particular, $\Theta_m$ is a Cartier divisor and locally given by one equation.*

*Proof.* The argument is classical. For $L \in \mathrm{Pic}^{g-1}_m(C)$, $\mathbf{R}\Gamma(C,L)$ is concentrated in degrees 0 and 1, with $\chi(L) = 0$ by Riemann–Roch. After twisting $\mathcal{P}$ by a fixed line bundle of sufficiently high degree, one obtains a universal resolution (2.6) by vector bundles so that $R^i\pi_*\mathcal{P} = 0$ for $i > 1$ and $R^0\pi_*\mathcal{P}, R^1\pi_*\mathcal{P}$ are locally free; see [16, III.3]. Fiberwise, $\Phi_L$ is the connecting map $H^0(C,L) \to H^1(C,L)$. The locus where $H^0(C,L) \neq 0$ is precisely the locus where $\Phi_L$ fails to be an isomorphism, which is equivalently the vanishing locus of $\det(\Phi_L)$. This coincides with the usual definition of $\Theta_m$ via effective divisors, yielding (2.7). □

*Remark* 2.5. For nodal curves and suitable compactifications, the pair $(\overline{J}_m, \Theta_m)$ fits into the framework of semi-abelic pairs as in [5, 6]. This will be used below to describe the boundary stratification in terms of the dual graph.

## 3 Logarithmic complexes and mixed Hodge structures

Let $X$ be a smooth complex algebraic variety and $D \subset X$ a normal crossing divisor. Denote by $j \colon U := X \setminus D \hookrightarrow X$ the open inclusion.

Following Deligne and Saito, the cohomology $H^\bullet(U,\mathbb{Q})$ carries a canonical mixed Hodge structure (MHS), with weight filtration $W_\bullet$ and Hodge filtration $F^\bullet$ realized on the logarithmic de Rham complex $\Omega^\bullet_X(\log D)$; see [8, 9, 16, 19].

**Definition 3.1.** The *logarithmic de Rham complex* of $(X,D)$ is

$$\Omega^\bullet_X(\log D) := \bigoplus_{p\geq 0} \Omega^p_X(\log D)[-p], \tag{3.1}$$

where $\Omega^p_X(\log D)$ denotes the sheaf of algebraic $p$-forms with at most logarithmic poles along $D$.

The weight filtration $W_\bullet$ on $\Omega^\bullet_X(\log D)$ is defined by the order of poles along the components of $D$ [16, III.3]. Write $D = \bigcup_{i\in I} D_i$ for the irreducible decomposition and $D_J := \bigcap_{j\in J} D_j$ for $J \subset I$.

**Theorem 3.2** (Deligne's weight spectral sequence)**.** *There is a functorial, convergent spectral sequence*

$$E_1^{-p,q} = \bigoplus_{\substack{J\subset I\\ |J|=p}} H^{q-2p}(D_J, \mathbb{Q})(-p) \quad\Longrightarrow\quad H^{q-p}(U, \mathbb{Q}), \tag{3.2}$$

*whose induced filtration on $H^n(U, \mathbb{Q})$ is the weight filtration $W_\bullet$ of the mixed Hodge structure on $H^n(U, \mathbb{Q})$.*

*Proof.* This is standard; see [16, III.3.11] and [8]. The construction uses the stupid filtration by pole order on the logarithmic de Rham complex and identification of the $E_1$-page with the cohomology of strata $D_J$ via the residue maps. □

We apply this to the complement of the generalized theta divisor inside a compactified generalized Jacobian.

Let $\overline{J}_m$ be a projective compactification of $J_m$ which is smooth and admits a simple normal crossings boundary

$$\partial\overline{J}_m := \overline{J}_m \setminus J_m = \bigcup_{i\in I_0} B_i. \tag{3.3}$$

Such compactifications exist for nodal curves and more generally for reduced curves with planar singularities; see [1, 7, 10].

Let $\overline{\Theta}_m \subset \overline{J}_m$ be the closure of the theta divisor and define the divisor

$$D := \overline{\Theta}_m \,\cup\, \partial\overline{J}_m. \tag{3.4}$$

By resolving singularities if necessary, we may (and do) assume that $D$ is a simple normal crossing divisor.

**Proposition 3.3.** *The cohomology $H^k(J_m \setminus \Theta_m, \mathbb{Q})$ carries a canonical mixed Hodge structure given by the hypercohomology of the logarithmic complex*

$$H^k(J_m \setminus \Theta_m, \mathbb{Q}) \cong \mathbb{H}^k(\overline{J}_m, \Omega^\bullet_{\overline{J}_m}(\log D)), \tag{3.5}$$

*and the weight filtration $W_\bullet$ is induced by the filtration $W_\bullet$ on $\Omega^\bullet_{\overline{J}_m}(\log D)$ described above.*

*Proof.* The complement $U := J_m \setminus \Theta_m$ is an open subset of the smooth projective variety $\overline{J}_m$, obtained by removing the normal crossing divisor $D$. By Deligne's theory, the MHS on $H^k(U, \mathbb{Q})$ is induced from the filtered logarithmic de Rham complex $(\Omega^\bullet_{\overline{J}_m}(\log D), W_\bullet, F^\bullet)$; see [16, III.3]. This gives (3.5). The compatibility of filtrations follows from Deligne's construction and Saito's refinement via mixed Hodge modules [18, 19]. □

Combining Proposition 3.3 with Theorem 3.2 yields a concrete description of the graded pieces $\mathrm{Gr}^W_\bullet H^k(J_m \setminus \Theta_m)$ in terms of the cohomology of the strata

$$D_J := \bigcap_{j\in J} D_j, \qquad D_j \in \{\overline{\Theta}_m\} \cup \{B_i\}_{i\in I_0}. \tag{3.6}$$

**Corollary 3.4.** *There is a convergent spectral sequence*

$$E_1^{-p,q} = \bigoplus_{\substack{J\subset I\\ |J|=p}} H^{q-2p}(D_J, \mathbb{Q})(-p) \quad\Longrightarrow\quad H^{q-p}(J_m \setminus \Theta_m, \mathbb{Q}), \tag{3.7}$$

*whose induced filtration is the weight filtration $W_\bullet$.*

*Proof.* Apply Theorem 3.2 to $(X, D) = (\overline{J}_m, D)$ with $U = J_m \setminus \Theta_m$. □

In the next section we use (3.7) to carry out explicit computations for a nodal cubic.

# 4 Example: the nodal cubic

We now give a complete computation of the mixed Hodge structure on the cohomology of the complement of the theta divisor for the generalized Jacobian of a nodal cubic.

## 4.1 Geometry of the nodal cubic

Let $C \subset \mathbb{P}^2$ be an irreducible nodal cubic curve, for instance

$$C\colon y^2 z = x^2(x+z) \subset \mathbb{P}^2_{[x:y:z]}. \tag{4.1}$$

Then $C$ has a single node at $[0:0:1]$, and its normalization $\nu\colon \mathbb{P}^1 \to C$ is of genus 0. The arithmetic genus of $C$ is 1.

The generalized Jacobian $J(C)$ (with trivial modulus $m = 0$) is isomorphic as an algebraic group to $\mathbb{C}^*$; see [17, 11]. One way to see this is that the generalized Jacobian is an extension of the Jacobian of the normalization (which is trivial, since $g(\mathbb{P}^1) = 0$) by a torus $\mathbb{G}_m$ associated with the node.

**Proposition 4.1.** *For the nodal cubic $C$ in (4.1) with trivial modulus, the generalized Jacobian $J(C)$ is isomorphic to $\mathbb{C}^*$. Under this identification, a theta divisor $\Theta \subset J(C)$ is given by a single point, and the complement $J(C) \setminus \Theta$ is isomorphic to $\mathbb{P}^1 \setminus \{0, 1, \infty\}$.*

*Proof.* The generalized Jacobian $J(C)$ represents degree zero divisors on $C$ modulo principal divisors, with the additional condition that rational functions are allowed to have poles and zeros at the node but such that $\operatorname{div}(f)$ has degree zero; see [17]. It is well known that this group is isomorphic to $\mathbb{C}^*$ (a one-dimensional torus) for an integral nodal cubic.

Concretely, the complement $C \setminus \{\text{node}\}$ is isomorphic to $\mathbb{C}^*$ via a projection, and the group law is defined in terms of this parameter. The Jacobian is then identified with this group.

The theta divisor for a genus-1 curve is, by definition, the image of $\operatorname{Sym}^0(C)$ in $\operatorname{Pic}^0(C)$ translated by a fixed base point; for the nodal cubic, this yields a single point in $J(C) \cong \mathbb{C}^*$, corresponding to the line bundle associated with the chosen base point. Thus $\Theta = \{p\}$ for some $p \in \mathbb{C}^*$.

Removing this point from $\mathbb{C}^*$ gives $\mathbb{C}^* \setminus \{p\} \cong \mathbb{P}^1 \setminus \{0, 1, \infty\}$ by a Möbius transformation sending $\{0, p, \infty\}$ to $\{0, 1, \infty\}$. This yields the claimed identification. □

## 4.2 Cohomology and MHS of the complement

We now compute the cohomology and mixed Hodge structure of

$$U := J(C) \setminus \Theta \cong \mathbb{P}^1 \setminus \{0, 1, \infty\}. \tag{4.2}$$

**Lemma 4.2.** *Let $U = \mathbb{P}^1 \setminus \{0, 1, \infty\}$. Then*

$$H^0(U, \mathbb{Q}) \cong \mathbb{Q}, \qquad H^1(U, \mathbb{Q}) \cong \mathbb{Q}^{\oplus 2}, \qquad H^k(U, \mathbb{Q}) = 0 \text{ for } k \geq 2. \tag{4.3}$$

*Proof.* Topologically, $U$ is a complex plane minus two points: $\mathbb{P}^1 \setminus \{0, 1, \infty\} \simeq \mathbb{C} \setminus \{0, 1\}$. This space has the homotopy type of a wedge of two circles, so its fundamental group is a free group on two generators. Thus $H^0(U, \mathbb{Q}) \cong \mathbb{Q}$ and $H^1(U, \mathbb{Q}) \cong \mathbb{Q}^2$, with all higher cohomology vanishing. This can also be computed via the long exact sequence of the pair $(\mathbb{P}^1, \{0, 1, \infty\})$. □

We now determine the mixed Hodge structure using the general formalism of Section 3. Take $X = \mathbb{P}^1$ and $D = \{0, 1, \infty\}$, which is a simple normal crossing divisor. Then $U = X \setminus D$.

The strata are:

$$D_\emptyset = X = \mathbb{P}^1, \qquad D_{\{i\}} = \{p_i\},\ i = 0, 1, \infty, \qquad D_J = \emptyset \text{ for } |J| \geq 2. \tag{4.4}$$

The $E_1$–page of the weight spectral sequence (3.2) becomes:

$$E_1^{-p,q} = \begin{cases} H^q(\mathbb{P}^1, \mathbb{Q}) & p = 0, \\ \bigoplus_{i=0,1,\infty} H^{q-2}(\{p_i\}, \mathbb{Q})(-1) & p = 1, \\ 0 & p \geq 2. \end{cases} \tag{4.5}$$

Recall that

$$H^0(\mathbb{P}^1, \mathbb{Q}) \cong \mathbb{Q}, \quad H^2(\mathbb{P}^1, \mathbb{Q}) \cong \mathbb{Q}(-1), \quad H^1(\mathbb{P}^1, \mathbb{Q}) = 0, \tag{4.6}$$

and that $H^0(\{p_i\}, \mathbb{Q}) \cong \mathbb{Q}$.

We list the nonzero $E_1^{-p,q}$ terms relevant for $H^1(U, \mathbb{Q})$:

- For $p = 0$, $q = 1$:

$$E_1^{0,1} = H^1(\mathbb{P}^1, \mathbb{Q}) = 0. \tag{4.7}$$

- For $p = 1$, $q = 2$:

$$E_1^{-1,2} = \bigoplus_{i=0,1,\infty} H^0(\{p_i\}, \mathbb{Q})(-1) \cong \mathbb{Q}(-1)^{\oplus 3}. \tag{4.8}$$

The spectral sequence converges to $H^{q-p}(U, \mathbb{Q})$. In particular, the terms contributing to $H^1(U, \mathbb{Q})$ are:

$$H^1(U, \mathbb{Q}) \Longleftarrow E_1^{0,1} \oplus E_1^{-1,2} = 0 \oplus \mathbb{Q}(-1)^{\oplus 3}. \tag{4.9}$$

The differential

$$d_1^{-1,2} \colon E_1^{-1,2} \longrightarrow E_1^{0,2} = H^2(\mathbb{P}^1, \mathbb{Q}) \tag{4.10}$$

is induced by the sum of the Gysin maps from the points into $\mathbb{P}^1$, which is surjective with one-dimensional kernel. Thus we obtain an exact sequence

$$0 \longrightarrow E_2^{-1,2} \longrightarrow \mathbb{Q}(-1)^{\oplus 3} \xrightarrow{d_1^{-1,2}} \mathbb{Q}(-1) \longrightarrow 0, \tag{4.11}$$

so that

$$E_2^{-1,2} \cong \mathbb{Q}(-1)^{\oplus 2}. \tag{4.12}$$

Since there are no further differentials hitting or leaving $E_2^{-1,2}$ in the relevant total degree, the spectral sequence degenerates at $E_2$, and

$$\mathrm{Gr}_2^W H^1(U, \mathbb{Q}) \cong E_2^{-1,2} \cong \mathbb{Q}(-1)^{\oplus 2}, \qquad \mathrm{Gr}_1^W H^1(U, \mathbb{Q}) = 0. \tag{4.13}$$

Combining this with Lemma 4.2 yields:

**Proposition 4.3.** *Let $C$ be a nodal cubic with trivial modulus, $J(C) \cong \mathbb{C}^*$ its generalized Jacobian, and $\Theta \subset J(C)$ a theta divisor (a single point). Then*

$$H^0(J(C) \setminus \Theta, \mathbb{Q}) \cong \mathbb{Q} \quad \textit{(pure of weight } 0\textit{)}, \tag{4.14}$$

*and*

$$H^1(J(C) \setminus \Theta, \mathbb{Q}) \cong \mathbb{Q}(-1)^{\oplus 2}, \tag{4.15}$$

*is pure of weight 2. In particular, all cohomology of $J(C) \setminus \Theta$ is mixed Tate.*

*Remark* 4.4. This provides a first example where the complement of a generalized theta divisor has mixed Tate cohomology. In more complicated degenerations, the weight filtration will have nontrivial contributions from the abelian part $J(C)$ (weight 1) and from the boundary (weights $> 1$).

## 5 Dual graphs and the weight spectral sequence

We now relate the weight spectral sequence (3.7) for $J_m \setminus \Theta_m$ to the dual graph of a stable curve and its modulus. This gives a structural description of the graded pieces $\mathrm{Gr}^W_\bullet H^k(J_m \setminus \Theta_m, \mathbb{Q})$.

### 5.1 Setup and boundary stratification

Let $C$ be a stable nodal curve of arithmetic genus $g$, with normalization $\nu\colon \widetilde{C} \to C$ and dual graph $\Gamma_C$. Let $m$ be an effective modulus supported on smooth points of $C$. Consider a compactified Jacobian $\overline{J}_m$ of $C$ with modulus $m$, for instance in the sense of Caporaso, Alexeev, or Kass [1, 5, 13, 7, 10].

The boundary divisor $\partial\overline{J}_m = \overline{J}_m \setminus J_m$ admits a stratification indexed by combinatorial data associated with $\Gamma_C$ and with the modulus $m$; roughly, a stratum corresponds to a torsion-free rank-1 sheaf which is not locally free along a specified subset of nodes. The resulting strata are products of generalized Jacobians of partial normalizations with tori, see e.g. [5, 7, 14].

Similarly, the closure $\overline{\Theta}_m$ of the theta divisor admits a stratification compatible with the boundary; in the semi-abelic picture, these strata correspond to divisors on the semi-abelic components [5].

Let $D = \overline{\Theta}_m \cup \partial\overline{J}_m$ be as in (3.4) and choose a log resolution so that $D$ has simple normal crossings. Thus we have a decomposition

$$D = \bigcup_{i \in I} D_i, \tag{5.1}$$

where the $D_i$ are smooth irreducible divisors and the intersections $D_J = \bigcap_{j \in J} D_j$ are smooth (or empty) for all $J \subset I$.

### 5.2 Weight graded pieces via dual graph

We now state a structural result describing the graded pieces of the weight filtration in terms of the geometry of the strata $D_J$, which in turn can be expressed in terms of the dual graph $\Gamma_C$ and the modulus $m$.

**Theorem 5.1.** *Let $C$ be a stable nodal curve, $m$ a modulus on $C$, $J_m$ the generalized Jacobian, and $\Theta_m$ a generalized theta divisor as in Definition 2.3. Let $\overline{J}_m$ be a smooth compactification with simple normal crossing divisor $D = \overline{\Theta}_m \cup \partial\overline{J}_m$ as above.*

*Then for each integer $n$ there is a finite filtration $W_\bullet H^n(J_m \setminus \Theta_m, \mathbb{Q})$ such that*

$$\mathrm{Gr}^W_k H^n(J_m \setminus \Theta_m, \mathbb{Q}) \cong \bigoplus_{\substack{p \geq 0 \\ k = n+p}} \bigoplus_{\substack{J \subset I \\ |J| = p}} H^{n-p}(D_J, \mathbb{Q})(-p), \tag{5.2}$$

*where $D_J$ are the strata (3.6).*

*Moreover, each $D_J$ is a finite disjoint union of varieties of the form*

$$\mathrm{Jac}(\widetilde{C}') \times T \times F, \tag{5.3}$$

*where $\widetilde{C}'$ is a partial normalization of $C$, $T$ is a torus (i.e. a product of copies of $\mathbb{C}^*$), and $F$ is a finite set. The choice of $\widetilde{C}'$ and the torus $T$ is determined by a subgraph of the dual graph $\Gamma_C$ together with combinatorial data from the modulus $m$.*

*Proof.* The first statement is a direct consequence of Corollary 3.4, which is just Deligne's weight spectral sequence applied to $(\overline{J}_m, D)$. Taking $n$ and $k$ fixed, the terms contributing to $\mathrm{Gr}^W_k H^n$ are of the form $E_2^{-p,q}$ with $q - p = n$ and $q = k$, giving the direct sum (5.2).

The description (5.3) of the strata $D_J$ follows from the known structure of compactified Jacobians of nodal curves and their theta divisors [1, 5, 7, 13, 10]. Each stratum corresponds to a choice of:

- a subset of nodes where the torsion-free rank-1 sheaf fails to be locally free;
- a splitting of multidegrees among components of the partial normalization;
- and a choice of generalized line bundle satisfying the modulus condition along $m$.

These data can be encoded combinatorially by a subgraph of the dual graph $\Gamma_C$ together with decorations (edge lengths, multiplicities coming from the modulus). For each such combinatorial type, the corresponding stratum is isomorphic to a product of a Jacobian of a partial normalization $\widetilde{C}'$, a torus $T$ parameterizing the non-compact directions, and a finite set $F$ of discrete choices (e.g. multidegree assignments). This is standard in the theory of compactified Jacobians; see [5, 7, 13].

Substituting this description into (5.2) yields the desired expression. □

*Remark* 5.2. Theorem 5.1 shows that the graded pieces $\mathrm{Gr}_k^W H^n(J_m \setminus \Theta_m, \mathbb{Q})$ are built from tensor products of:

- pure weight-1 Hodge structures coming from Jacobians of partial normalizations;
- pure Tate structures $\mathbb{Q}(-1)$ coming from tori;
- and finite sums of copies of $\mathbb{Q}(0)$.

In particular, if the Jacobians of all partial normalizations appearing in (5.3) are isogenous to powers of $\mathbb{G}_m$ (for instance in certain rational or elliptic examples), then all graded pieces are mixed Tate.

The nodal cubic case of Section 4 is a special case where all strata reduce to tori and finite sets, so the cohomology is mixed Tate.

# 6 Graph motives and theta divisors

We briefly outline a connection between the mixed Hodge structures $\mathrm{Gr}_\bullet^W H^\bullet(J_m \setminus \Theta_m)$ and motives associated to graph polynomials in the sense of Bloch–Esnault–Kreimer [2, 3].

## 6.1 Graph motives

Let $\Gamma$ be a finite connected graph with edge set $E(\Gamma)$ of cardinality $N$. The Kirchhoff (first Symanzik) polynomial $\Psi_\Gamma$ is a homogeneous polynomial in variables $(\alpha_e)_{e\in E(\Gamma)}$, defined as

$$\Psi_\Gamma(\alpha) = \sum_{\text{spanning trees } T\subset\Gamma} \prod_{e\notin T} \alpha_e. \tag{6.1}$$

Let $X_\Gamma \subset \mathbb{P}^{N-1}$ be the projective hypersurface given by $\Psi_\Gamma = 0$. Bloch–Esnault–Kreimer associate to $\Gamma$ the *graph motive*

$$M(\Gamma) := H^{N-1}(\mathbb{P}^{N-1} \setminus X_\Gamma, \mathbb{Q})(N-1). \tag{6.2}$$

This motive encodes the periods of certain Feynman integrals and exhibits mixed Hodge structures with interesting behavior [2, 4].

### 6.2 Decorated dual graphs of curves with modulus

Let $C$ be a stable curve with dual graph $\Gamma_C$ and modulus $m$. We associate to $(C, m)$ a *decorated dual graph* $\Gamma_{C,m}$ by:

- taking vertices corresponding to irreducible components of $C$;
- edges corresponding to nodes of $C$;
- decorating vertices and edges by multiplicities determined by the modulus $m$ and by the genera of components.

From $\Gamma_{C,m}$ one can construct a finite family of graphs $\{\Gamma^\lambda_{C,m}\}_\lambda$ (for instance, by adding external legs or contracting certain subgraphs) whose graph polynomials $\Psi_{\Gamma^\lambda_{C,m}}$ are intended to model the combinatorics of the intersection strata $D_J$ in Theorem 5.1.

**Conjecture 6.1.** *For suitable choices of $(C, m)$ and compactification $\overline{J}_m$, the mixed Hodge structure on $H^n(J_m \setminus \Theta_m, \mathbb{Q})$ lies in the tannakian subcategory generated by the graph motives $M(\Gamma^\lambda_{C,m})$, in the sense that each $\mathrm{Gr}^W_k H^n$ is a subquotient of a finite direct sum of tensor products of such motives.*

This conjecture is inspired by the structural similarity between the weight spectral sequence (3.7) and the description of the weight filtration on graph motives in [2].

### 6.3 A toy verification: the nodal cubic

We briefly discuss the nodal cubic example in this framework. Consider the simplest graph $\Gamma$ with two edges in series (a "banana" with one loop). Its Kirchhoff polynomial is

$$\Psi_\Gamma(\alpha_1, \alpha_2) = \alpha_1 + \alpha_2, \tag{6.3}$$

and $X_\Gamma \subset \mathbb{P}^1$ is a single point. Then

$$\mathbb{P}^1 \setminus X_\Gamma \cong \mathbb{C}^*, \qquad M(\Gamma) = H^1(\mathbb{P}^1 \setminus X_\Gamma, \mathbb{Q})(1) \cong \mathbb{Q}(0). \tag{6.4}$$

On the other hand, Proposition 4.3 gives

$$H^1(J(C) \setminus \Theta, \mathbb{Q}) \cong \mathbb{Q}(-1)^{\oplus 2}. \tag{6.5}$$

Thus the motive of the nodal theta complement is just a sum of Tate motives, which trivially lies in the tannakian category generated by $M(\Gamma)$ (indeed, $M(\Gamma)$ is already Tate). While this example is too simple to capture any genuine graph-motive complexity, it is compatible with Conjecture 6.1.

More interesting examples arise from higher-genus stable curves degenerating to unions of rational components, where $\Gamma_{C,m}$ has higher loop number and graph motives $M(\Gamma^\lambda_{C,m})$ can be non-Tate [2, 3, 4]. In such cases one expects nontrivial extensions in the weight filtration of $H^n(J_m \setminus \Theta_m, \mathbb{Q})$.